\documentclass[11pt]{article}

\usepackage[utf8x]{inputenc}

\usepackage[font=small]{caption}
\usepackage{amsmath,amsthm,amssymb,stmaryrd,bigints,mathabx}

\usepackage[pdftex]{graphicx}
\usepackage{subcaption}
\usepackage[titletoc,title]{appendix}
\usepackage{datetime,color,cancel,listings,authblk}
\usepackage[hyperfootnotes=false]{hyperref}
\hypersetup{
  colorlinks = true,
  linkcolor = black,
  citecolor = blue,
  urlcolor = black
}
\usepackage[margin=2.0cm]{geometry}

\usepackage{float}

\floatstyle{ruled}
\newfloat{algorithm}{thp}{lop}
\floatname{algorithm}{Algortithm}

\newcommand{\beq}{\begin{equation}}
\newcommand{\ee}{\end{equation}}
\newcommand{\bac}{\begin{array}{c}}
\newcommand{\ea}{\end{array}}
\newcommand{\bal}{\begin{aligned}}
\newcommand{\eal}{\end{aligned}}

\newcommand{\real}{\operatorname{Re}}
\newcommand{\imag}{\operatorname{Im}}

\newcommand{\vertiii}[1]{{\left\vert\kern-0.25ex\left\vert\kern-0.25ex\left\vert #1 
    \right\vert\kern-0.25ex\right\vert\kern-0.25ex\right\vert}}

\begin{document}
\newtheorem{theorem}{Theorem}[section]
\newtheorem{lemma}[theorem]{Lemma}
\newtheorem{remark}[theorem]{Remark}
\newtheorem{observation}[theorem]{Observation}
\newtheorem{definition}[theorem]{Definition}
\newtheorem{example}[theorem]{Example}
\newtheorem{corollary}[theorem]{Corollary}
\newtheorem{assumption}{Assumption}
\newtheorem{property}{Property}


\title{Modulation instability and  convergence of the random phase approximation for stochastic sea states}

\author{Agissilaos G. Athanassoulis\thanks{Corresponding author email: aathanassoulis@dundee.ac.uk} }
\author{Irene Kyza}

\affil{Department of Mathematics, University of Dundee}

\maketitle

\begin{abstract} The nonlinear Schr\"odinger equation is widely used as an approximate model for the evolution in time of the water wave envelope. In the context of simulating ocean waves, initial conditions are typically generated from a measured power spectrum using the random phase approximation, and periodized on an interval of length $L$. It is known that most realistic ocean waves power spectra do not exhibit modulation instability, but the  most severe ones do; it is thus a natural question to ask whether the periodized random phase approximation has the correct stability properties. In this work we specify a random phase approximation scaling so that, in the limit of $L\to\infty,$ the stability properties of the periodized problem are identical to those of the continuous power spectrum on the infinite line. Moreover, it is seen through concrete examples that using a too short computational domain can completely suppress the modulation instability.
\end{abstract}

\medskip 



\section{Introduction}

Ocean wave models can be classified as {\em phase-resolved} or {\em phase-averaged}. A phase-resolved model aims to describe each individual wave crest and trough -- for example the first principles, fully nonlinear, free-boundary problem for potential flow, or asymptotic approximations thereof (Boussinesq, KdV, Zakharov, nonlinear Schr\"odinger etc) are all phase-resolved models. On the other hand, phase-averaged models (such the CSY, Alber and Hasselmann  equations) aim to describe a {\em stochastic sea state} on a  macroscopic and/or coarse-grained level. In that sense any phase-averaged model is a statistical model of waves, and typically uses some additional assumptions  in order to come up with closed equations for a phase-averaged representation. Phase-resolved models can be very accurate in terms of hydrodynamics. 
However, they still have to be initialized, and any prediction with them ultimately is as good as the initial conditions used. Since the vast majority of {\em data} available for real-world ocean waves is phase-averaged  (typically power spectra / autocorrelations), the crisp exactness of ideal hydrodynamics is out of reach anyway when modelling the ocean. On the other hand, {\em phase-averaged models have the advantage of being designed  with real-world data in mind}. They are well suited for the large space and time scales required in ocean engineering, and are known to provide valuable insights for the qualitative behaviour of stochastic wavefields. Phase-resolved simulations can also be used in conjunction with data to study stochastic sea states in a Monte Carlo sense, i.e. when used with a large population of initial conditions ({\em realizations} of the sea state) generated from phase-averaged measurements. This is more accurate,  but also much more computationally expensive than phase-averaged models -- and thus reserved for validation or cases of special interest. Typically, a measured power spectrum (itself a phase-averaged representation of the sea surface elevation) is used to generate a population of initial conditions on a finite computational domain of size $L$  through the {\em random-phase approximation}.

In this work we study the convergence for the random phase approximation in terms of qualitative behaviour, i.e. whether it correctly captures the presence or not of modulation instability in the original problem. We work in a regime where it is known that both modulation instability and its absence are possible \cite{Alber1978,Andrade2020b,Athanassoulis2018,Dysthe2008,Gramstad2017,Ribal2013a}. 


\subsection{Phase-averaged models and the modulation instability}

In 1962, Nobel-prize winner K. Hasselmann published what is now known as the Hasselmann equation,  a {phase-averaged}  equation describing the energy transfer between various wavenumbers \cite{Hasselmann1962a}. 
 The starting point of \cite{Hasselmann1962a} is a stochastic sea state, where both the background and perturbations are {\em gaussian and homogeneous random processes}. Both gaussianity and homogeneity are widely considered to be fundamental features of deep-water gravity waves \cite{Ochi1998}; this is largely consistent with empirical observation, and is furthermore supported intuitively by an informal application of the Central Limit Theorem.  Today many state of the art approaches still use the Hasselmann equation, or similar moment equations (often collectively called  ``wave kinetic  equations'', WKE) to study different aspects of real-world ocean waves, including long-distance propagation \cite{wam1988,ericrogers2020,Janssen2008},  extreme events \cite{Fedele2015,Fedele_2010_pof,fnsr009} etc. Moreover, it has inspired a number of rigorous mathematical studies not necessarily closely tied to real-world ocean waves. These include the well-posedness theory of generalized wave kinetic  equations \cite{germain2020}, wave turbulence \cite{Buckmaster2021,skipp23,el2023}, as well as the rigorous investigation of the closure conditions involved in the derivation of the WKE \cite{ZHani23,Deng2021,Deng2021a}.

A few years after the introduction of Hasselmann's equation, what is now called {\em modulation instability (MI)} entered the picture,  in the study of the nonlinear Schr\"odinger equation (NLS) as an approximate model for ocean waves \cite{zakharov1968stability,Benjamin1967} (groups working on different problems were also encountering  aspects of the MI around the same time, see \cite{Zakharov2009} for a historical overview). The presence of MI in focusing dispersive waves means that plane waves are inherently unstable. It is now known that the MI is a pretty universal feature of focusing, dispersive waves \cite{Biondini2016,Biondini2018}. Moreover, \cite{Biondini2016,Biondini2018} point out that the  MI can be seen as {\em the tendency of  inhomogeneities to initially grow}, and subsequently settle in a particular pattern.   It is natural to ask whether stochastic  wavefields with sufficiently narrow power spectra behave like plane waves in that respect. The fact that ocean waves are in many cases narrowband  makes this a very natural question.

A first attempt to investigate how the MI might show up in a phase-averaged moment equation (analogous to Hasselmann's equation) was made by Longuet-Higgins \cite{LH1976}. In that work Davey-Stewartson dynamics were used and it was found that  the MI for a plane wave  would indeed shows up in the moment equations. The investigation  for a wavefield with a narrow but realistic power spectrum (not a delta-function)  was first performed in \cite{Alber1978}. A sufficient condition for  exponential growth of small inhomogeneities was derived, in the form of a nonlinear system of equations depending on the power spectrum (what I. E. Alber called ``the integral stability eigenvalue relation'' \cite{Alber1978}). We will call this kind of instability of a stochastic sea state  {\em generalized modulation instability (gMI)}, as opposed to the {\em classical MI (MI)} which is the growth of inhomogeneities around a  plane wave only. It turns out that the MI can be seen as a special case of the gMI for a delta-function power spectrum (which practically amounts to a plane wave) \cite{Athanassoulis2017}. The notable feature of Alber's analysis is that sufficiently narrow and intense spectra can be seen to be {\em unstable} (meaning that inhomogeneities would spontaneously grow), while sufficiently broad and low-intensity spectra are expected to be {\em stable} (meaning that inhomogeneities do not grow; in fact later it was proved that inhomogeneities disperse in the stable case \cite{Athanassoulis2018}). The quantitative resolution of the stability condition proved to require the development of  novel numerical and analytical techniques, and took several decades to mature \cite{Andrade2020b,Athanassoulis2018,Ribal2013a,Gramstad2017} -- more on that below.

The possibility  of realistic sea states for which gMI is present (i.e. homogeneity is not stable) has profound consequences in the modelling of ocean waves. This was first articulated in the seminal paper \cite{Crawford1980}, where a second moment equation is derived for water waves with an assumption of gaussianity only -- not homogeneity. This is called the CSY equation, and it is a broadband, 2-dimensional two-space moment equation, derived from the Zakharov equation using only a gaussian closure. It is the most accurate moment equation for water waves in the literature.  Assuming homogeneity is a robust feature of the sea state, then one can go to sixth order effects for the homogeneous wavefield, and essentially reproduce Hasselmann's equation \cite{Crawford1980}. When a narrowband approximation is used in \cite{Crawford1980}, it is found that the question of whether homogeneity is preserved can be investigated quantitatively, and Alber's condition is essentially recovered. Crucially, if homogeneity is not stable, then the effects of the gMI are second order.   In \cite{Andrade2020b} this point was further strengthened:  assuming that the initial sea state in the CSY equation consists of  an $O(1)$ homogeneous part and an $o(1)$ inhomogeneous part, a CSY stability condition was derived, controlling whether the inhomogeneous part grows exponentially or not. This was done without any narrowband approximation, and it yields the broadband, 2-dimensional counterpart of Alber's ``stability relation''.

So \cite{Alber1978,Crawford1980,Andrade2020b} collectively present a fundamental dichotomy: starting with the a priori assumptions of gaussianity and quasi-homogeneity of a given sea state (meaning the initial presence of a small inhomogeneity),  an Alber-type stability condition for that particular sea state can be produced. If it is found to be stable, then homogeneity is robust, the inhomogeneity is going to disperse, and one can safely go to higher-order effects in order to derive homogeneous kinetic wave equations - as in, e.g., \cite{Hasselmann1962a,Janssen2003}. In that sense the Alber stability of the sea state is {\bf a compatibility condition}, allowing the self-consistent use of homogeneity.

\begin{figure}
\subfloat{%
  \includegraphics[clip,width=0.4\textwidth]{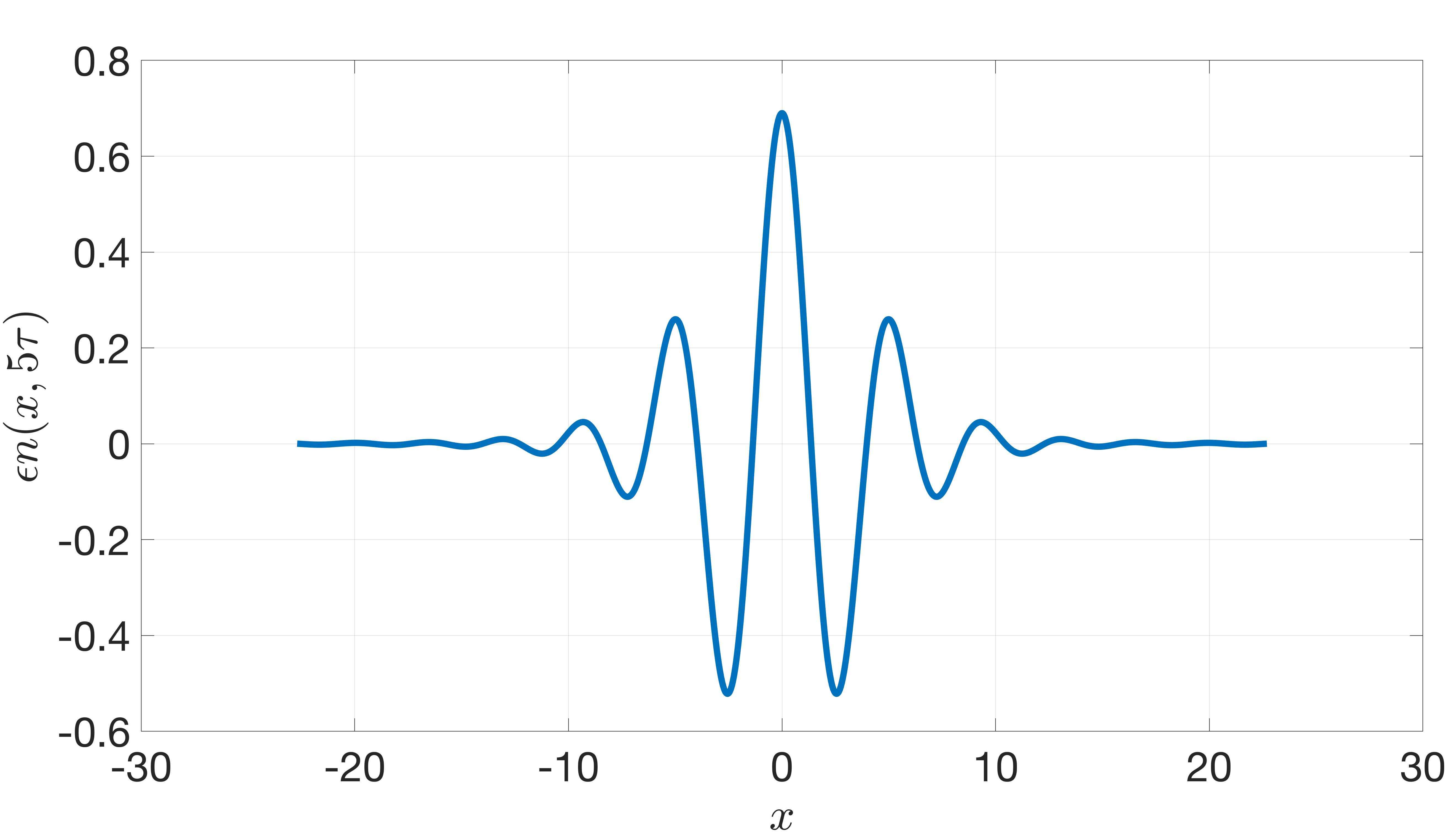}%
}\,\,
\subfloat{%
  \includegraphics[clip,width=0.57\textwidth]{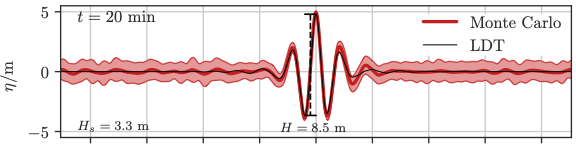}%
  \vspace{0.44cm}
}
\caption{{\bf Left:} Inhomogeneity produced by the unstable modes of the Alber equation for a JONSWAP spectrum \cite{Athanassoulis2017}. {\bf Right:} Average profile of extreme events for a JONSWAP spectrum; Monte Carlo simulation using modified NLS dynamics (image provided by T. Grafke) \cite{Dematteis2017}.}
\label{fig:new1}
\end{figure}

If the spectrum is found to be unstable however, then any inhomogeneity (which in the field could be due to wave breaking, gusts of wind, localized extreme events, presence of ships etc) will be spontaneously amplified, and this is a leading order effect. Thus homogeneity cannot be consistently used as a permanent structural feature of an unstable sea state: {\em the homogeneous sea state is basically an unstable equilibrium.} These  sea states where gMI is present have to be studied with different methods, and are expected to behave very differently from the stable ones.
This pumping of energy in inhomogeneities for unstable spectra has been reproduced in many works, including \cite{Dysthe2008,dtk2003,GramstadTrulsen2007,Onorato2006Exper,moriyasuda,Athanassoulis2023,Onorato2003,Ribal2013a}. In the presence of gMI, large localized extreme events  emerge which clearly break homogeneity.  It is worth noting that working out the fundamental scalings of the gMI can give a qualitative description of these localized extreme events; see Figure \ref{fig:new1} for a comparison of asymptotic analysis of the Alber equation \cite{Athanassoulis2017} with state of the art Monte Carlo results \cite{Dematteis2017}.

The remaining question is the quantitative resolution of the stability condition: when is a given power spectrum expected to be stable or unstable? 
The papers \cite{Andrade2020b,Athanassoulis2018,Ribal2013a,Gramstad2017} all investigate the onset of gMI for JONSWAP parametric spectra. The relevant parametrization here is in $\alpha,$ a parameter directly affecting the severity of the sea state, and $\gamma,$ a parameter primarily affecting how narrow the spectrum is. In practice, more severe sea states (storms) also tend to have narrower spectra -- i.e. both $\alpha$ and $\gamma$ tend to increase or decrease together.  While different techniques (and in some cases different equations) are used, the findings are broadly in agreement: the vast majority of realistic ocean wave spectra are stable, but the most extreme sea states on record  cross over to the unstable region. On one hand, this is consistent with the widespread success of homogeneous approaches and the absence of MI-like effects in mild and moderate seas \cite{fnsr009,Lee2022}. On the other hand, it is consistent with large, dangerous rare events  requiring a separate investigation  \cite{Olagnon,Ochi1998,Kharif2003a}. 
It is worth noting that investigation of non-parametric spectra \cite{Athanassoulis2021} is also possible, and yields similar results. 

\subsection{Onset of  instability and Monte Carlo setup}
\label{sec:realiz}

Understanding these unstable sea states is crucial for several reasons: on one hand, it is likely that they play a key role in real-world oceanic rogue waves; on the other hand they represent a mathematical regime that is very poorly understood and in fact most of the time we simply {\em assume them away}. Because of the instability of homogeneity, any phase-averaged model that assumes permanent homogeneity cannot be used as a first-principles model. For that we should use Monte Carlo simulations of the underlying phase-resolved problem. 

That is, we have to use phase-averaged data (power-spectra) to create a large number of realizations of the sea surface elevation. These initial realizations will be evolved in time over a finite computational (or experimental) domain of length $L$. Very often periodic boundary conditions are used (although different kinds of boundary conditions also exist in the literature, e.g. non-reflective absorbing etc).  It is standard practice to do this with the random-phase approximation \cite{Ochi1998,Komen1994}, basically a superposition of plane waves with random phases and amplitudes consistent with the power spectrum. 

Let us assume nonlinear Schr\"odinger  (NLS) dynamics governing the evolution in time of a complex wave envelope $u(x,t).$ The original physical problem is typically framed in infinite space \cite{Chiang2005},
\begin{equation}\label{eq:nls1}
	i \partial_t U + \frac{p}2 \Delta U + \frac{q}2 |U|^2U=0, \quad x\in \mathbb{R}, \quad  \mathop{\sup}\limits_{x\in\mathbb{R}}|U(x,t)|< \infty,
\end{equation}
with initial data $U_0(x)$ being a realization of a homogeneous random process with power spectrum $S(k).$
On the other hand, the truncated problem that we can readily simulate is the periodized NLS
\begin{equation}\label{eq:nls2a}
\begin{array}{c}
i \partial_t u + \frac{p}2 \Delta u + \frac{q}2 |u|^2u=0, \qquad
x\in [-\frac{L}2,\frac{L}2], \\
u(-\frac{L}2,t)=u(\frac{L}2,t),  \qquad 
	\partial_x u(-\frac{L}2,t)=\partial_x u(\frac{L}2,t),
\end{array}
\end{equation}
with initial data of the form $u_0(x)=\sum A_j e^{2\pi i (k_jx+\phi_j)}$ where the $\phi_j$ are i.i.d. random variables uniformly distributed in $[0,1).$

The main question that we address in this paper is the following:
\begin{itemize}
	\item Under what conditions can we guarantee that the presence or not of gMI in the original problem \eqref{eq:nls1} is accurately captured in the periodized problem \eqref{eq:nls2a}?
\end{itemize}
We answer this question by deriving a periodized Alber equation for the truncated problem, along with its stability relation which contains explicitly  the lengthscale $L.$ Then we can see that the instability condition for the truncated problem converges to that of the original problem if appropriate scalings are used (roughly speaking as long as $L$ is large enough and $O(L)$ plane waves are used -- see Theorem \ref{thrm:main}). More specifically, we find a way to create realizations of the sea state that leads to the correct stability condition. 
To highlight the role  of $L$ in the qualitative behaviour of the solution, we investigate its role in the fully nonlinear MI in Section \ref{sec:2}.

Recent theoretical breakthroughs suggest that different preparation of the initial phase-resolved sea state may lead to completely different regimes \cite{ZHani23,el2023}. Thus, the preparation of the initial data deserves more detailed attention. This is already underway on the engineering side, where  initialization and calibration is now starting to be seen  as an iterative workflow  in the context of high quality data \cite{penesis23} rather than a standard formula. 

\section{Main results}

\begin{theorem}[Periodized Alber equation] \label{thrm:main} Let $S(k)$ be  a smooth power spectrum
with $\mathop{supp} S(k) \subseteq [0,k_{max}],$ and consider
 equation \eqref{eq:nls2a} with stochastic initial data $u(x,0)=u_0(x)$ where
\begin{equation}\label{eq:scaling}
u_0(x) = \sum\limits_{j=1}^{M} A_j e^{2\pi i (\frac{jx}L + \phi_j) }, \qquad A_j = \sqrt{ \frac{1}L S(\frac{j}L) } , \qquad \phi_j \sim \mbox{ i.i.d. } U(0,1), \qquad M=k_{max}\cdot L.
\end{equation}
The Alber equation for the evolution in time of inhomogeneity $\rho(x,y,t)$ of this problem is
\[
\bac
	i\partial_t \rho + \frac{p}2 \left(\Delta_x - \Delta_y\right) \rho
+q  \left[\Gamma(x-y)+\rho(x,y,t)\right]\left[ \rho(x,x,t) - \rho(y,y,t) \right]  =0, \quad  x,y \in [\frac{L}2,\frac{L}2],
\ea
\]
equipped with periodic boundary conditions,
where
\[
\Gamma(x) = \sum\limits_{n} \frac{1}L S(\frac{n}L) e^{2\pi i \frac{xn}L},
\]
and the linear
instability  condition for this problem is
\begin{equation}\label{eq:tony}
\exists \xi\in\mathbb{Z} \qquad \exists \omega \in \{ z\in \mathbb{C}: \real(z)>0\} \qquad\qquad \widetilde{h}_{L}(\xi,\omega)=1
\end{equation}
where 
\begin{equation}\label{eq:thehscal}
\widetilde{h}_{L}(\xi,\omega)=\frac{iq}L  \sum\limits_{k\in \mathbb{Z}}  \frac{    S(\frac{k-\xi}{L}) -  S(\frac{k}{L})  }{ \omega  +ip\frac{2\pi^2 \xi(2k-\xi)}{L^2} }.
\end{equation}
Moreover, in the limit $L\to\infty,$ $\xi\to\infty,$ $\xi/L\to X$
\begin{equation}\label{eq:convmain}
\lim\limits_{L\to\infty} h_L(LX,\omega)= 
iq \int\limits_{k=-\infty}^{+\infty}
\frac{ 
 S(k- \frac{X}{2}) -  S(k+\frac{X}{2}) 
}{
  \omega   +4\pi^2 i p k X 
  } dk.
\end{equation}
\end{theorem}

\begin{corollary}[Convergence of the instability condition as $L\to\infty$]
Under the assumptions of Theorem \ref{thrm:main},  the linear instability condition \eqref{eq:tony} of the periodized problem converges  to the linear instability condition for the infinite line problem, \cite{Alber1978,Athanassoulis2018},
\begin{equation}\label{eq:ilsc}
	\exists X\in\mathbb{R} \qquad \exists \omega \in \{ z\in \mathbb{C}: \real(z)>0\} \qquad\qquad iq \int\limits_{k=-\infty}^{+\infty}
\frac{ 
 S(k- \frac{X}{2}) -  S(k+\frac{X}{2}) 
}{
  \omega   +4\pi^2 i p k X 
  } dk=1,
\end{equation}
as $L\to\infty.$
\end{corollary}

\begin{remark}
We don't claim that equation \eqref{eq:scaling} is the only way to recover the correct stability condition. Other random-phase-type strategies to create realizations of the sea state are also used in the literature, and it seems likely that they can also be calibrated lead to the same stability condition.  
\end{remark}

%
%
%
%
%
%
%
%
%
%
%

The rest of the paper is organized as follows: 
Sections \ref{sec:3.1} and \ref{sec:3.2} go over the derivation of the periodized Alber equation and its stability relation. They are written in an accessible way that helps motivate and frame the proof, which can be found in Section \ref{sec:3.3}. Methods to check the stability condition are recalled in Section \ref{sec:real4}. We explore the role of $L$ in the MI and gMI in Section \ref{sec:2} - in particular, it is seen that  a short computational domain can completely stabilize the MI as well as the gMI.
Ramifications  are discussed in Section \ref{sec:conc}.

\section{Stability analysis for periodized stochastic sea states} \label{sec:4}

\subsection{Derivation of the stability condition}\label{sec:3.1}

Let us briefly recall the derivation of the Alber equation in the infinite space case.
Assuming NLS dynamics for the wave envelope on infinite space \eqref{eq:nls1}
and following \cite{Alber1978,Athanassoulis2018}, we can derive the second moment equation
\begin{equation}\label{eq:010}
	i\partial_t R + \frac{p}2 \left(\Delta_x - \Delta_y\right) R
+q R(x,y,t) \left[ R(x,x,t) - R(y,y,t) \right]  =0
\end{equation}
 for the two-space autocorrelation, 
\begin{equation}
	R=R(x,y,t)=E[u(x,t)\overline{u}(y,t)].
\end{equation}
The  standard gaussian closure
\[
E[u(x, t) \bar{u}(x, t) u(x, t) \bar{u}(y, t)]=2 E[u(x, t) \bar{u}(y, t)] E[u(x, t) \bar{u}(x, t)]
\]
is used to arrive to equation \eqref{eq:010}, cf. Appendix B of \cite{Athanassoulis2018} for more details. Moreover, by separating the  wavefield  into a homogeneous and non-homogeneous part 
\begin{equation}
R(x,y,t)= \Gamma(x-y)+\rho(x,y,t), 
\end{equation}
we derive the standard (infinite space) Alber equation for the evolution of the inhomogeneous part,
\begin{equation}\label{eq:01}
	i\partial_t \rho + \frac{p}2 \left(\Delta_x - \Delta_y\right) \rho
+q  \left[\Gamma(x-y)+\rho(x,y,t)\right]\left[ \rho(x,x,t) - \rho(y,y,t) \right]  =0.
\end{equation}
Observe that if $\rho=0,$ equation \eqref{eq:01} is automatically satisfied. In other words, any homogeneous autocorrelation $R(x,y,t)=\Gamma(x-y)$ is automatically a solution of equation \eqref{eq:010}; the point is to investigate if a particular homogeneous solution  is  stable under inhomogeneous perturbations (in which case an initially small inhomogeneity  $\rho$ is guaranteed to stay small and even disperses) or unstable  (in which case any inhomogeneity is expected to grow rapidly).

Now let us consider the periodization of the NLS  on a torus of length $L.$
Starting from equation \eqref{eq:nls2a} and performing the same steps as before, the same Alber equation \eqref{eq:01} can still be derived  --  now equipped with periodic boundary conditions instead of non-growth at infinity.
Moreover, the periodicity of $u$ is inherited by its second moments,
and thus we can express moments in terms of Fourier series on the interval $[-\frac{L}2,\frac{L}2],$
\begin{equation}\label{eq:fcoeffs}
	\Gamma(y) = \sum\limits_{n\in\mathbb{Z}} P_n e^{2\pi i \frac{ny}L}, \qquad 	\rho(x,y,t) = \sum\limits_{k,l \in \mathbb{Z}} r_{k,l}(t) e^{2\pi i \frac{kx+ly}L}.
\end{equation}
For now we treat the Fourier coefficients $P_n$ as known. We will revisit their relationship to the continuous power spectrum $S(k)$ in the next subsection.
Now assuming quasi-homogeneity, i.e.
\begin{equation}
	\rho(x,y,0)=o(1),
\end{equation}
we linearise around homogeneity
(i.e. drop  $O(\rho^2)$ terms)  and express everything in terms of the  Fourier coefficients introduced in \eqref{eq:fcoeffs}:
\begin{equation}\label{eq:p2}
\begin{array}{c}
i \sum\limits_{k,l \in \mathbb{Z}} \partial_tr_{k,l}(t) e^{2\pi i \frac{kx+ly}L} + \frac{p}2 \sum\limits_{k,l \in \mathbb{Z}} r_{k,l}(t) \left((2\pi i \frac{k}L)^2 - (2\pi i \frac{l}L)^2\right) e^{2\pi i \frac{kx+ly}L}+ \\
+q  \left(\sum\limits_{n \in \mathbb{Z}} P_n e^{2\pi i \frac{n(x-y)}L}\right)
\left(  \sum\limits_{k,l \in \mathbb{Z}} r_{k,l}(t) \left[e^{2\pi i \frac{k+l}L x} - e^{2\pi i \frac{k+l}L y} \right]\right)  =0. 
\end{array}
\end{equation}

By taking the inner product with $e^{2\pi i \frac{k'x+l'y}L}$ (and suppressing the primes) we obtain
%
%
%
%
\begin{equation}\label{eq:poikjuyh}
  \partial_tr_{k,l}(t)  +i \frac{2\pi^2 p}{L^2}  r_{k,l}(t) (k-l)(k+l) 
-iq \left[  P_{-l} -  P_{k} \right]  \sum\limits_{K\in \mathbb{Z}}   r_{K,k+l-K}(t) =0.
\end{equation}
Denote for brevity 
\begin{equation}
\vec r(t) := \{ r_{k,l}(t)\}_{k,l\in\mathbb{Z}},
\end{equation}
now equation \eqref{eq:poikjuyh} is linear and can be summarised as
\begin{equation}\label{erq:A}
	\frac{d}{dt} \vec r(t) =  A \vec r(t),
\end{equation}
in other words 
 the growth of $|r(t)|$ in time can be inferred by the eigenvalues of the matrix $A.$  This can provide a powerful black-box method (see \cite{Andrade2020b} for a similar approach in the CSY equation), but would not allow us, e.g., to show convergence as  $L\to\infty.$ Thus, we will proceed to work with the Laplace transform as in \cite{Athanassoulis2018}.

To help understand better what affects the eigenvalues, first of all we bring equation \eqref{eq:poikjuyh} to mild form
\begin{equation}\label{eq:p04}
\begin{array}{c}
r_{k,l}(t) = e^{-ip\frac{2\pi^2(k+l)(k-l)}{L^2} t} r_{k,l}(0) +  \int\limits_{\tau=0}^t e^{-ip\frac{2\pi^2(k+l)(k-l)}{L^2} (t-\tau)} iq \left[  P_{-l} -  P_{k} \right]  \sum\limits_{K\in \mathbb{Z}}   r_{K,k+l-K}(\tau) d\tau ,
\end{array}
\end{equation}
and proceed to define
\begin{equation}
f(\xi,t)=\sum\limits_{K\in \mathbb{Z}}   r_{K,\xi-K}(t), \qquad \xi\in\mathbb{Z};
\end{equation}
this corresponds to the $\xi$'th Fourier coefficients of the inhomogeneous part of the position density $|u(x,t)|^2$. 
(The point is that  $f(\xi,t)$ really drives the problem, in the sense that if we now $f(\xi,t)$ then $r_{k,l}(t)$ can be computed by substitution in equation \eqref{eq:p04}.) 

Now set 
\[
k+l=\xi \implies l=\xi-k, \qquad k-l=2k-\xi;
\]
thus, equation  \eqref{eq:p04} becomes
\begin{equation}\label{eq:pp05}
r_{k,\xi-k}(t) =
 e^{-ip\frac{2\pi^2 \xi (2k-\xi)}{L^2} t} r_{k,\xi-k}(0) + 
  \int\limits_{\tau=0}^t e^{-ip\frac{2\pi^2\xi(2k-\xi)}{L^2} (t-\tau)} 
  iq \left[  P_{k-\xi} -  P_{k} \right]  \sum\limits_{K\in \mathbb{Z}}   r_{K,\xi-K}(\tau) d\tau.
  \end{equation}
By summing equation \eqref{eq:pp05} in $k$ it follows that
\begin{equation}\label{eq:LT1}
f(\xi,t) =
\phi_L(\xi,t) + 
  \int\limits_{\tau=0}^t h_L(\xi,t-\tau)  f(\xi,\tau) d\tau,
\end{equation}
where
\begin{equation}
	\phi_L(\xi,t)=\sum\limits_k e^{-ip\frac{2\pi^2 \xi (2k-\xi)}{L^2} t} r_{k,\xi-k}(0), \qquad 
	h_L(\xi,t-\tau) = \sum\limits_k e^{-ip\frac{2\pi^2\xi(2k-\xi)}{L^2} (t-\tau)} 
  iq \left(  P_{k-\xi} -  P_{k} \right).
\end{equation}
This is a simplified form of \eqref{eq:poikjuyh} in the sense that we have one discrete variable $\xi$ in the new unknown function $f$ instead of two, $k,l$ in the old variable $\vec r$, while we still can recover the full $\vec r(t)$ through
\begin{equation}\label{eq:prewaveop}
r_{k,l}(t) = 
 e^{-ip\frac{2\pi^2(k+l)(k-l)}{L^2} t} r_{k,l}(0) +  \int\limits_{\tau=0}^t e^{-ip\frac{2\pi^2(k+l)(k-l)}{L^2} (t-\tau)} iq \left[  P_{-l} -  P_{k} \right]  f(k+l,\tau) d\tau
\end{equation}
(in particular, if $f(t)$ does not grow exponentially, neither does $\vec r(t)$).

By Laplace transform equation \eqref{eq:LT1}  becomes
\begin{equation}
\widetilde{f}(\xi,\omega) = \widetilde{\phi}_L(\xi,\omega)+ \widetilde{h}_L(\xi,\omega)\widetilde{f}(\xi,\omega) \implies \widetilde{f}(\xi,\omega) = \frac{\widetilde{\phi}_L(\xi,\omega)}{1-\widetilde{h}_L(\xi,\omega)}.
\end{equation}
Thus, the existence of poles in the right half plane,
\[
\exists \omega^* \,\, : \real(\omega^*)>0 \,\, \mbox{ and } \widetilde{h}_L(\xi,\omega^*)=1
\]
means that  $f(\xi,t)$ exhibits exponential growth.
In terms of working out explicitly the function $\widetilde{h}_L,$ observe that for any $\omega$ with $\real \omega >0$ we have
\begin{equation}\label{eq:001}
\begin{array}{c}
\widetilde{h}_L(\xi,\omega) =\int\limits_{t=0}^\infty e^{-\omega t} \left[ \sum\limits_{k\in\mathbb{Z}} e^{-ip\frac{2\pi^2\xi(2k-\xi)}{L^2} t} 
  iq \left(  P_{k-\xi} -  P_{k} \right) \right] dt=
iq \sum\limits_{k\in\mathbb{Z}}  \frac{    P_{k-\xi} -  P_{k} }{ \omega +ip\frac{2\pi^2\xi(2k-\xi)}{L^2} } 
\end{array}
\end{equation}
and similarly 
\begin{equation}\label{eq:002}
\widetilde{\phi}_L(\xi,\omega)	= \sum\limits_{k\in\mathbb{Z}} \frac{r_{k,\xi-k}(0)}{\omega + ip\frac{2\pi^2\xi(2k-\xi)}{L^2}},
\end{equation}
so that finally
\begin{equation}\label{eq:003}
\bac
\widetilde{f}(\xi,\omega)	=\frac{
\sum\limits_{k\in\mathbb{Z}} \frac{r_{k,\xi-k}(0)}{\omega + ip\frac{2\pi^2\xi(2k-\xi)}{L^2}}
}{
1- iq \sum\limits_{k\in\mathbb{Z}}  \frac{    P_{k-\xi} -  P_{k} }{ \omega +ip\frac{2\pi^2\xi(2k-\xi)}{L^2} }
 }.
\ea
\end{equation}

\subsection{Discretizing continuous power spectra}\label{sec:3.2}

In the previous subsection we took the Fourier coefficients $P_n$ for granted. Here we  consider how they are related to the original continuous power spectrum $S(k).$
Given a spectrum $S(k)$ a standard (periodic) realization of the  envelope corresponding to it can be generated through
\begin{equation}\label{eq:ansatzrpa}
u_0(x) = \sum\limits_{j=1}^{M} A_j e^{2\pi i ({k_jx} + \phi_j) }, \qquad A_j = \sqrt{ \,\delta k\, S(k_j) } , \qquad \phi_j \sim U(0,1), \qquad k_j = j \cdot \delta k,
\end{equation}
so that $A_j^2\approx  \int\limits_{k_{j-1}}^{k_j} S(k)dk$  \cite{Ochi1998}. 
Observe that in general there is some flexibility in the wavenumber discretization, 
\begin{equation}
 k_j = j\cdot \delta k, \quad j=1,2,\dots,M, \qquad \delta k = \frac{m}{L}, \qquad m,M \in \mathbb{N}.
\end{equation}
Now the autocorrelation of the random-phase process \eqref{eq:ansatzrpa} is
\[
\begin{array}{c}
\Gamma(x-y)=E[u(x)\overline{u}(y)] = \sum\limits_{j,j'=1}^M E[ A_j \bar{A}_{j'} e^{2\pi i[k_jx-k_{j'}y+\phi_j-\phi_{j'}]}] = 
\sum\limits_{j,j'=1}^M A_j \bar{A}_{j'} e^{2\pi i[k_jx-k_{j'}y]} E[  e^{2\pi i (\phi_j-\phi_{j'}  ) } ] = \\

= \sum\limits_{j,j'=1}^M A_j \bar{A}_{j'} e^{2\pi i[k_jx-k_{j'}y]} \delta_{j,j'} = \sum\limits_{j=1}^M |A_j|^2 e^{2\pi i k_j(x-y)}.
\end{array}
\]
Therefore
\begin{equation}\label{eq:Pn00}
\begin{array}{c}
P_n = \frac{1}L \int\limits_{y=-\frac{L}2}^{\frac{L}2} e^{-2\pi i \frac{ny}L} \sum\limits_{j=1}^M |A_j|^2 e^{2\pi i k_jy} dy = \frac{1}L \sum\limits_{j} |A_j|^2 \int\limits_{y=-\frac{L}2}^{\frac{L}2} e^{2\pi i (k_j-\frac{n}L)y} dy =  \sum\limits_{j} |A_j|^2 \delta_{k_j,\frac{n}L} =
\\
 =  \sum\limits_{j} \frac{m}L P(j\delta k) \delta_{j\delta k,\frac{n}L} = 
\left\{ \begin{array}{cl}
\frac{m}L  P(\frac{n}L), & n=mj \mbox{ for some } j \\
0, & n \mod m \neq 0
\end{array}
\right.
\end{array}
\end{equation}

Now let us go back to the parameters $m,M:$
the parameter $M$ controls how many distinct frequencies/wavenumbers are used, while $m$ controls how close these wavenumbers are taken. Thus the largest wavenumber included  is 
\[
k_{max} = M\cdot \delta k = \frac{Mm}{L}.
\]
Since we are interested in gravity waves it is reasonable to select a fixed maximum wavenumber, e.g. $k_{max}\approx10.$ Thus if $k_{max}=O(1)$ and $L\gg 1,$ this leads to
\begin{equation}
	M\cdot m = O(L),
\end{equation}
which in principle could allow for different wavenumber discretization strategies. Now if we return to the choice
\begin{equation}\label{eq:mscaling}
m=1, \qquad M=k_{max} \cdot L,
\end{equation}
that was made in the statement of Theorem \ref{thrm:main},
equation \eqref{eq:Pn00} implies that 
\begin{equation}\label{eq:Pn1}
\begin{array}{c}
P_n = \frac{1}L S(\frac{n}L).
\end{array}
\end{equation}

\subsection{Proof of  theorem \ref{thrm:main}} \label{sec:3.3}

In section \ref{sec:3.1} the expression for $\widetilde{h}_L(\xi,\omega)$ was derived in terms of the Fourier coefficients $P_n,$ cf. equation \eqref{eq:001}.
In section \ref{sec:3.2} it was shown how the scaling of equation \eqref{eq:scaling} for the initial data leads to the Fourier coefficients $P_n,$ cf. equation \eqref{eq:Pn1}. Combining equations \eqref{eq:001} and \eqref{eq:Pn1} leads to equation \eqref{eq:thehscal}.

Now, observe that
with the change of variables $k'=2k - \xi$ (which means that $k' \in  2\mathbb{Z} - \xi,$ $k=\frac{k'+\xi}2,$ $k-\xi=\frac{k'-\xi}2,$) we have
\begin{equation} \label{eq:resc2}
\bac
\displaystyle 
\widetilde{h}_L(\xi,\omega) = \frac{iq}L  \sum\limits_{k\in \mathbb{Z}}  \frac{    S(\frac{k-\xi}{L}) -  S(\frac{k}{L})  }{ \omega  +ip\frac{k2\pi^2 \xi(2k-\xi)}{L^2} }=

\frac{iq}2  \sum\limits_{k'\in 2\mathbb{Z}-\xi}  \frac{    S(\frac{k'-\xi}{2L}) -  S(\frac{k'+\xi}{2L})   }{ \omega   +4\pi^2 i p \frac{k'}{2L} \frac{\xi}L } \frac{2}{L}.
\ea
\end{equation}
Now denote $X\in\mathbb{R}$ a fixed number independent from $L$
and take the joint limit
\[
L\to\infty, \quad \xi\to\infty \qquad \mbox{ so that } \quad \frac{\xi}{L} \to X.
\]
In that regime $\widetilde{h}_L(\xi,\omega)$
can be seen as a Riemann sum converging  to
\begin{equation}\label{eq:36}
\bac
\lim \widetilde{h}_L = \frac{iq}2 \int\limits_{s=-\infty}^{+\infty}
\frac{ 
 S(\frac{s}2- \frac{X}{2}) -  S(\frac{s}2+\frac{X}{2}) 
}{
\omega   +4\pi^2 i p \frac{s}{2} X 
  } ds.
\ea
\end{equation}
Where $k'/L\to s$ and $2/L \to ds.$ This double-spacing
is because  $k'$ in equation \eqref{eq:resc2} spans the odd or even numbers only, depending on $\xi;$ that's why $2/L \to ds.$ Equation \eqref{eq:convmain} follows  with the change of variables $k=s/2$ in equation \eqref{eq:36}. 

With regard to the convergence of the Riemann integral in equation \eqref{eq:36}, note that as along as $\real(\omega)\neq 0,$ the denominator does not become zero, and the integrand inherits the smoothness and decay properties of the power spectrum $S(k).$
\qed

\section{Resolving the stability condition} \label{sec:real4}

The instability conditions both for the periodized problem and the infinite line problem (equation \eqref{eq:tony} and \eqref{eq:ilsc} respectively) are of the form
\begin{equation}\label{eq:AlbER2}
\exists \xi \qquad \exists \omega \in \{ z\in \mathbb{C}: \real(z)>0\} \qquad\qquad \widetilde{h}(\xi,\omega)=1
\end{equation}
for some known function $\widetilde{h}(\xi,\omega)$ which is by construction analytic in $\omega$ as long as $\real(\omega)>0.$ This is a system of two equations (the real and imaginary part of $\widetilde{h}$) in three unknowns ($\xi,$ $\real(\omega),$ $\imag(\omega)$) and as such determining the existence or not of solutions is not trivial.

A technique that can be used to determine the existence or not of solutions is based on the argument principle, and was analyzed in detail in \cite{Athanassoulis2018}. In a nutshell, it is based on the idea that, given a closed subset $D \subseteq \{ z\in\mathbb{C}: \real(z)>0\}$ with piecewise smooth boundary $\Gamma=\partial D \subseteq \{ z\in\mathbb{C}: \real(z)>0\},$  the curve $\widetilde{h}(\xi,\Gamma) = \{ \widetilde{h}(\xi,s) \in \mathbb{C}, \,\, s\in \Gamma\}$ circumscribes $\widetilde{h}(\xi,D)=\{ \widetilde{h}(\xi,s) \in \mathbb{C}, \,\, s\in D\}.$ 
A natural choice for $D$ is a right semicircle of radius $R\frac{1}\epsilon$ with the strip $\{ \real(Z)<\epsilon\}$ excluded. As $\epsilon \to 0,$ this set exhausts the open right half plane. So now we 
have $\Gamma$ being a D-shaped contour consisting of the straight line $\{ \epsilon + is\}_{s=-\frac{1}\epsilon}^{\frac{1}\epsilon}$ and the arc $\{\frac{1}\varepsilon(cos\theta + i \sin\theta)\}_{\theta=\mathop{arccos} \varepsilon^2}^{-\mathop{arccos} \epsilon^2}.$ 
Computing $\widetilde{h}(\xi,s)$ for $s\in \Gamma$ is straightforward; moreover, for small enough epsilon the values of $\widetilde{h}$ on the arc are going to be very close to zero, hence the nontrivial part is the images of the straight line -- this is analogous to Nyquist plots in control theory. 

A more elementary technique is just plotting all the values attained on a fine enough grid in the domain $D.$ Finally, when the existence of an unstable mode is detected, plotting the distance $|\widetilde{h}-1|$ can help locate its approximate location and work locally if necessary. Indeed, a number of numerical choices are required (numerical tolerances etc), it is often helpful to use several different methods in conjunction.

\section{MI and gMI can be completely suppressed on a short interval} 
\label{sec:2}

In this section we go over in some more detail the role of $L,$ the length of the computational domain, in the MI and gMI. We start by recalling in detail the MI for the real line in Section \ref{sec:MIR}, and for the interval of length $L$ in section \ref{sec:MIL}. Then we investigate numerically the nonlinear phase of the MI in Section \ref{sec:3}. We plot the inhomogeneity (in this case the difference from the plane wave background) and discuss its qualitative behaviour for various values of $L.$ In Section \ref{sec:nlgmit} we investigate numerically the nonlinear evolution of a small inhomogeneity added to a homogeneous background solution which serves as a simple model realization of a sea state with a very narrow spectrum. By plotting the inhomogeneity we see a picture completely analogous to the MI: a short computational domain $L$ can completely stabilize the problem, i.e. make the inhomogeneity stay small -- while on larger computational domains the inhomogeneity grows by two orders of magnitude. In all cases the inhomogeneity is numerically zero at the endpoints of the interval, while the background solution satisfies periodic boundary conditions.

\begin{figure}
\includegraphics[width=0.49\textwidth]{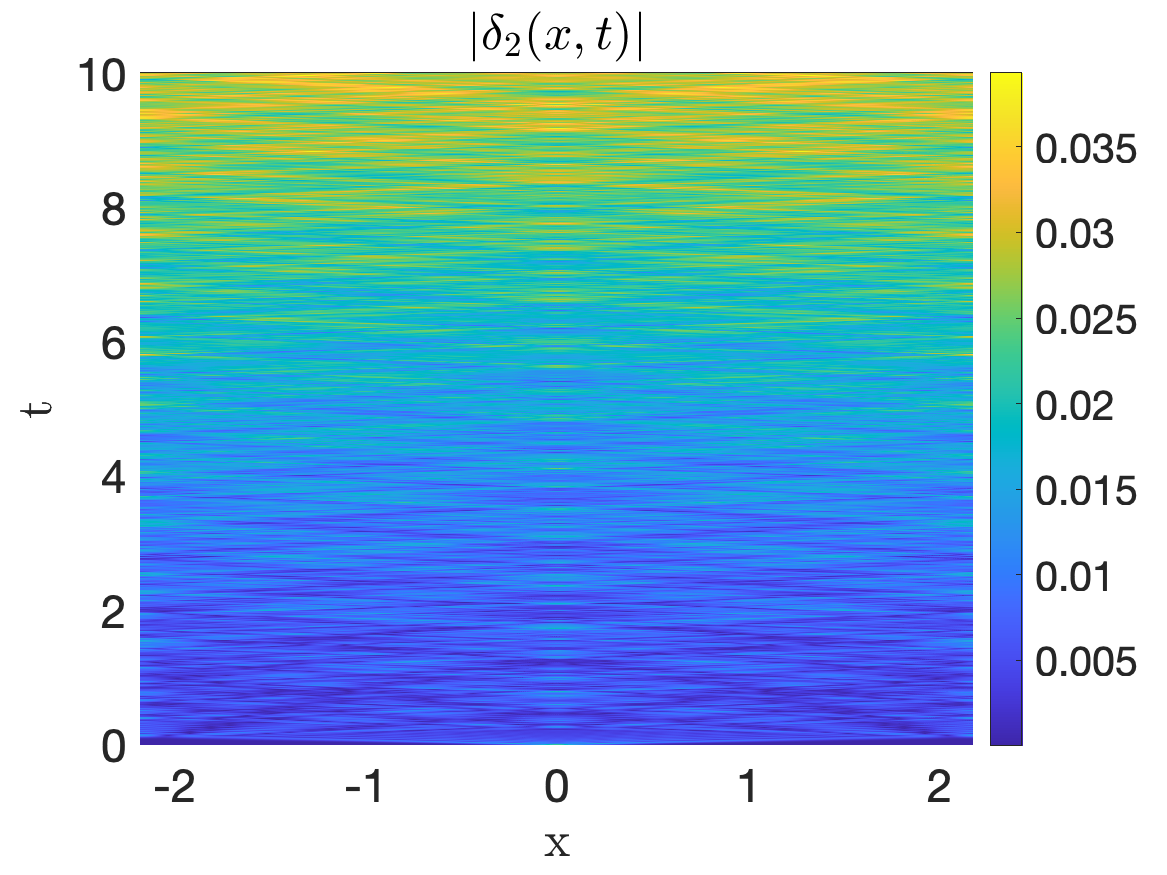} \,
\includegraphics[width=0.49\textwidth]{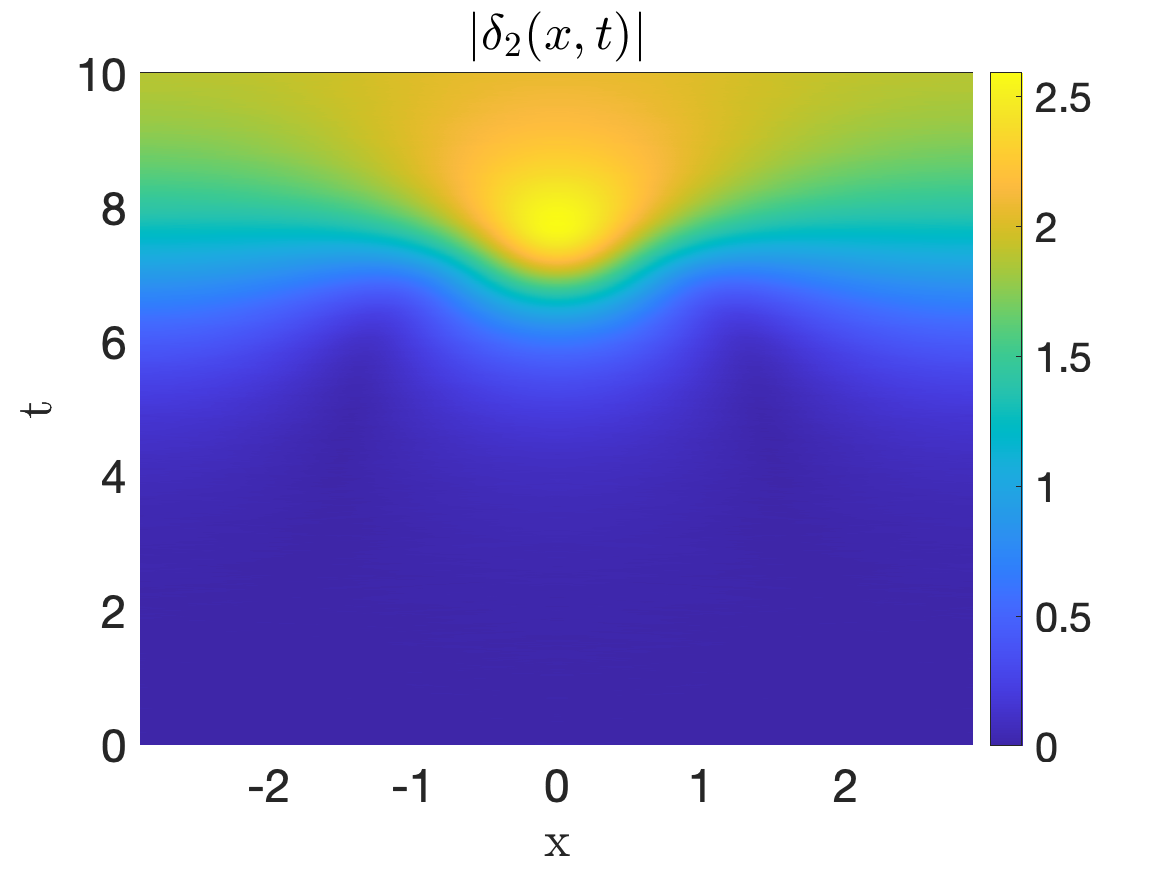} \\
\includegraphics[width=0.49\textwidth]{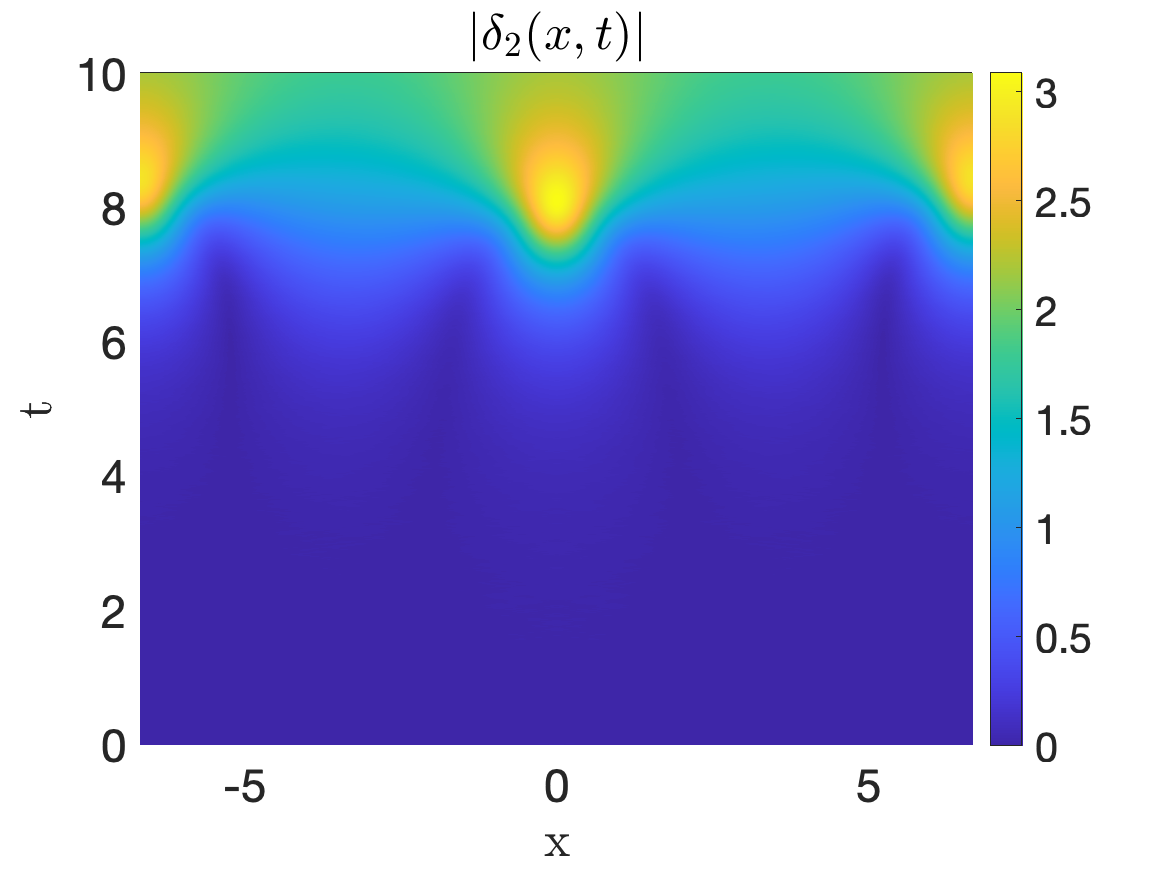} \,
\includegraphics[width=0.49\textwidth]{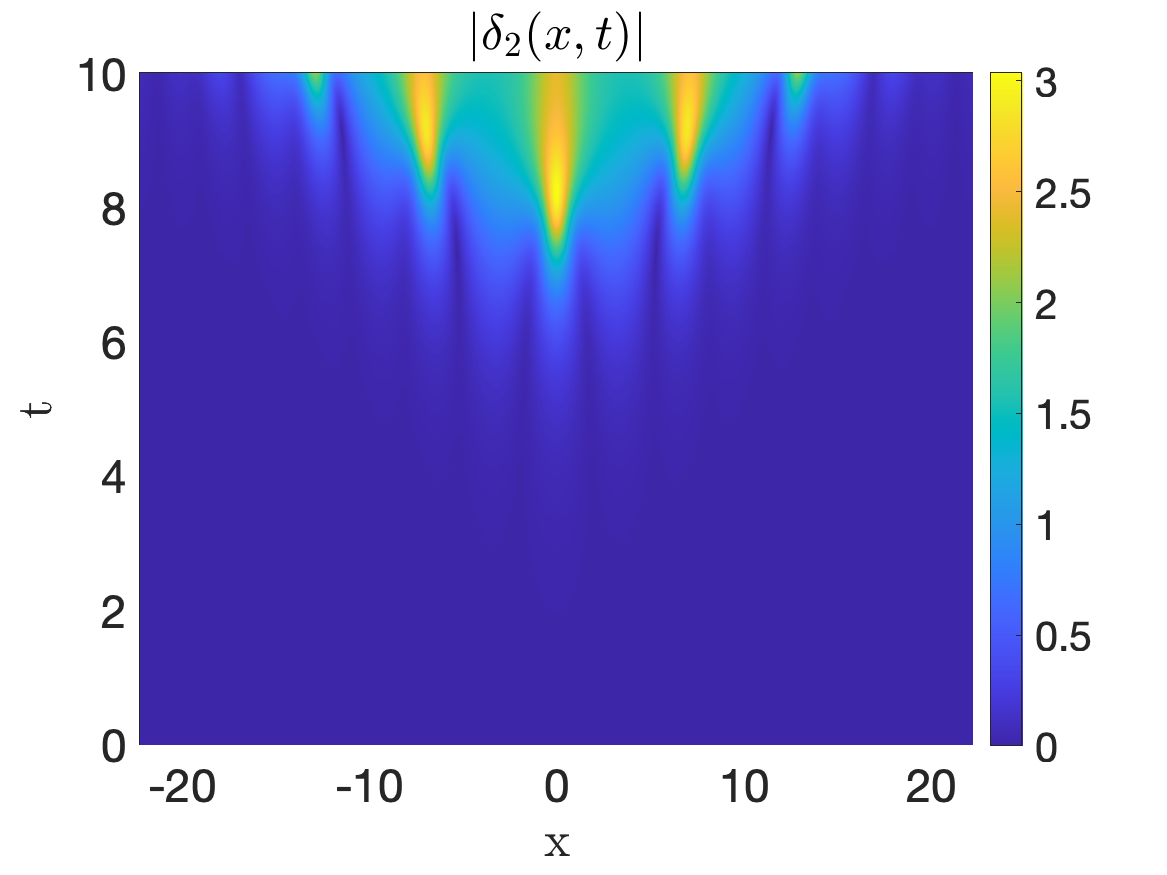}
\caption{Space-time plot of the modulus of the inhomogeneity, $|\delta_2(x,t)|,$ when propagated on intervals of different length $L$ according to equations \eqref{eq:nls4} -- \eqref{eq:defddd}. Growth of the inhomogeneity to $O(1)$ is evidence of MI. The bifurcation length $L_c$ is reported in equation \eqref{eq:L1}.  {\bf Top left:} $L=0.98L_c.$ {\bf Top right:} $L=1.3L_c.$ {\bf Bottom left:} $L=3L_c.$ {\bf Bottom right:} $L=10L_c.$}
\label{fig:b}
\end{figure}


\begin{table}
\begin{center}
\begin{tabular}{ |c|cccc| } 
 \hline
 & $N=0.98$ & $N=1.3$ & $N=3$ & $N=10$ \\
 \hline
$j=1$ & $0.0359$ & $2.59$ & $3.08$ & $3.03$ \\
$j=2$ & $0.0393$ & $2.59$ & $3.09$ & $3.03$ \\
$j=3$ & $0.149$ & $2.65$ & $3.17$ & $3.12$ \\
$j=4$ & $0.25$ & $2.69$ & $3.22$ & $3.18$ \\
$j=5$ & $0.03$ & $2.57$ & $3.05$ & $2.96$ \\
\hline
\end{tabular}
\end{center}
\caption{Maximum value of the inhomogeneity 
$\mathop{\mathrm{max}}\limits_{x,t} |\delta(x,t)|$ for $x\in [-L/2,L/2]$ and $t\in [0,10],$ with $L=NL_c$ (cf. equation \eqref{eq:L1}) and initial inhomogeneity $\delta_j(x,0)$ (cf. equation \eqref{eq:ICinh}).  } \label{tab:1}
\end{table}

\subsection{MI on the real line}
\label{sec:MIR}
We briefly go through the standard MI calculation on the real line, for equation \eqref{eq:nls1}. This will prepare the ground for the subsequent analysis on a bounded interval. 
Equation \eqref{eq:nls1} is known to be  well-posed in Zhidkov spaces \cite{Gallo2004,Zhidkov1987}; this is the natural framework for this discussion. Also, 
since we are dealing with the focusing NLS, we will assume without loss of generality that $p,q>0.$
Moreover, equation \eqref{eq:nls1} admits the exact solution
\begin{equation}\label{eqpwdef}
	w(x,t):= A e^{iqA^2 t},
\end{equation}
which is the simplest plane wave solution with amplitude $A>0$.
To study the stability of this plane wave solution one can consider whether small perturbations grow. This leads to the perturbed initial value problem
\begin{equation}\label{eq:nls2}
\begin{array}{c}
	i \partial_t u + p \Delta u + q |u|^2u=0,  \qquad 
u(x,0) = A(1+\delta_0(x))
\end{array}
\end{equation}
for some initial perturbation $\delta_0(x)$ which is small in an appropriate sense, $\delta_0=o(1)$. Using the substitution
\begin{equation}\label{eq:4}
	u(x,t) = A e^{iqA^2 t}(1+\delta(x,t)) 
\end{equation}
one readily computes that problem \eqref{eq:nls2} is equivalent to
\begin{equation}\label{eq:nls3}
\begin{array}{c}
i\delta_t + p\Delta \delta + q A^2 (\delta + \bar\delta) + q A^2 (\delta + \bar\delta) \delta + q A^2 |\delta|^2 (1+\delta),
\qquad 
\delta(x,0) = \delta_0(x).
\end{array}
\end{equation}
Dropping higher order terms, 
 we obtain the linearized problem for the perturbation, namely
\begin{equation}\label{eq:linearised}
\begin{array}{c}
i\delta_t + p\Delta \delta + q A^2 (\delta + \bar\delta),
\qquad 
\delta(x,0) = \delta_0(x).
\end{array}
\end{equation}
By expanding equation \eqref{eq:linearised} into its real and imaginary parts, and denoting 
\begin{equation}\label{eq:rimparts}
\delta(x,t) = \alpha(x,t) + i \beta(x,t),
\end{equation}
we eventually obtain the system
\begin{equation}\label{eq;sol0}
\begin{array}{c}
\partial_{tt}\beta + \left( p^2 \Delta\Delta +2pq A^2\Delta \right) \beta=0, \qquad 
\partial_t \alpha + p\Delta \beta=0.
\end{array}
\end{equation}
This now can be solved explicitly with separation of variables, leading to the construction of the modes
\begin{equation}\label{eq:solbetaR}
\begin{array}{c}
\zeta \in \mathbb{R}, \qquad
\beta_\zeta(x,t) = e^{i [\zeta x + \omega(\zeta) t]} + c.c., \qquad
 \omega^2(\zeta) = \zeta^2 [p^2 \zeta^2 - 2pq A^2].
\end{array}
\end{equation}
(Here c.c. stands for complex conjugate.)
More general solutions can be formed by superpositions of these modes, 
$\beta(x,t) = \int\limits M(\zeta) \beta_{\zeta}(x,t) d\zeta.$ 
The instability is due to the fact that, 
\begin{equation}\label{eq:cond}
\begin{array}{l}
|\zeta|<A\sqrt{2\frac{q}p}  \implies 
\omega(\zeta) =\pm i |\omega(\zeta)|
\end{array}
\end{equation}
and thus the corresponding modes $\beta_\zeta(x,t)$ of equation \eqref{eq:solbetaR}  contain an exponentially growing component.
Therefore the solution $\delta$ of the linearized equation \eqref{eq:linearised} generally grows exponentially in time, i.e. the plane wave solution of the NLS is linearly unstable. Moreover the unstable wavenumbers and their rate of growth follow from this analysis.

In this  problem there are three parameters, $p,q,A.$ However by rescaling the problem according to
$\tau = qA^2 \, t,$ $\chi = A\sqrt{\frac{q}p} \, x,$ $U(\chi,\tau) = \frac{1}A \, u(x,t),$
the equation is mapped to the ``canonical'' NLS
$i \partial_{\tau} U +  \partial_{\chi\chi} U +  |U|^2U=0,$
which has the plane wave solution $W(\chi,\tau)=e^{i\tau}.$ That is, the exact values of $p,q,A$ don't play an important role and the nature of the modulation instability is the same for any $p,q,A>0.$ This widely understood universality of the MI may have contributed to  an implicit expectation that the periodized problem automatically enjoys similar properties. 

\subsection{MI on an interval of length $L$}\label{sec:MIL}

Now let us consider the NLS equation on an interval of length $L$ equipped with periodic boundary conditions, namely \eqref{eq:nls2a}.
The steps described above in equations \eqref{eqpwdef}-\eqref{eq;sol0} were carried out with boundary conditions of boundedness at infinity; however, one readily checks that each step is equally valid with the periodic boundary conditions on $[-L/2,L/2].$ Indeed, the plane wave is still a solution, the algebra leading to \eqref{eq:nls3} is the same, and finally the linearization and separation of real and imaginary parts are the same. So now we have to solve the problem 
\begin{equation}\label{eq:s7856rgh}
\begin{array}{c}
\partial_{tt}\beta + \left( p^2 \Delta\Delta +2pq A^2\Delta \right) \beta=0, \quad 
\beta(-\frac{L}2,t)=\beta(\frac{L}2,t), \quad 
\partial_x\beta(-\frac{L}2,t)=\partial_x\beta(\frac{L}2,t).
\end{array}
\end{equation}
This is the first time that the parameter $L$ really comes into play, and separation of variables now leads to the {\em discrete modes}
\begin{equation}
\begin{array}{c}
\displaystyle 
	\beta_n(x,t) = e^{2\pi i(\frac{ n x}L +\omega_n t)} + c.c., \qquad  (2\pi \omega_n)^2 = (\frac{2\pi n}L)^2 [p^2(\frac{2\pi n}L)^2 -2pqA^2].
\end{array}
\end{equation}
While this is analogous to equation \eqref{eq:solbetaR}, there is a very important difference: {\em depending on the values of $p,q,A,L,$ there may not be even one unstable mode}. Indeed, $\omega_0=0,$ and for any $0\neq n \in \mathbb{Z}$ we have 
\begin{equation}
\omega_n^2<0 
\iff L> \frac{2\pi |n|}A \sqrt{\frac{p}{2q}}.
\end{equation}
That is, the computational domain $L$ has to be larger than
\begin{equation}\label{eq:L1}
L_c:=\frac{2\pi \sqrt{p}}{ A \sqrt{{2q}} },	
\end{equation}
otherwise  the MI is {\em completely suppressed}. 
Observe moreover that the lengthscale $L_c$
 does not depend on the initial inhomogeneity $\delta_0(x),$ and becomes larger when the  nonlinearity becomes weaker (i.e. when $q,A>0$ decrease). 
 For water waves with typical wavenumber $k_0$ (equivalently typical wavelength $\lambda_0=2\pi/k_0$), the coefficients are $p=\sqrt{g}/(8k_0^{3/2}),$ $q=\sqrt{g}k_0^{5/2}/2$ \cite{Chiang2005} leading to 
\begin{equation}\label{eq:l11}
L_c =  \frac{1}{4\pi\sqrt{2}} \frac{\lambda_0^2}{A}.
\end{equation} 
Note that the waves  that would be dangerous for ships in the ocean and  carry most of the surface wave energy have wavelengths in the hundreds of meters.

\begin{figure}
\includegraphics[width=0.49\textwidth]{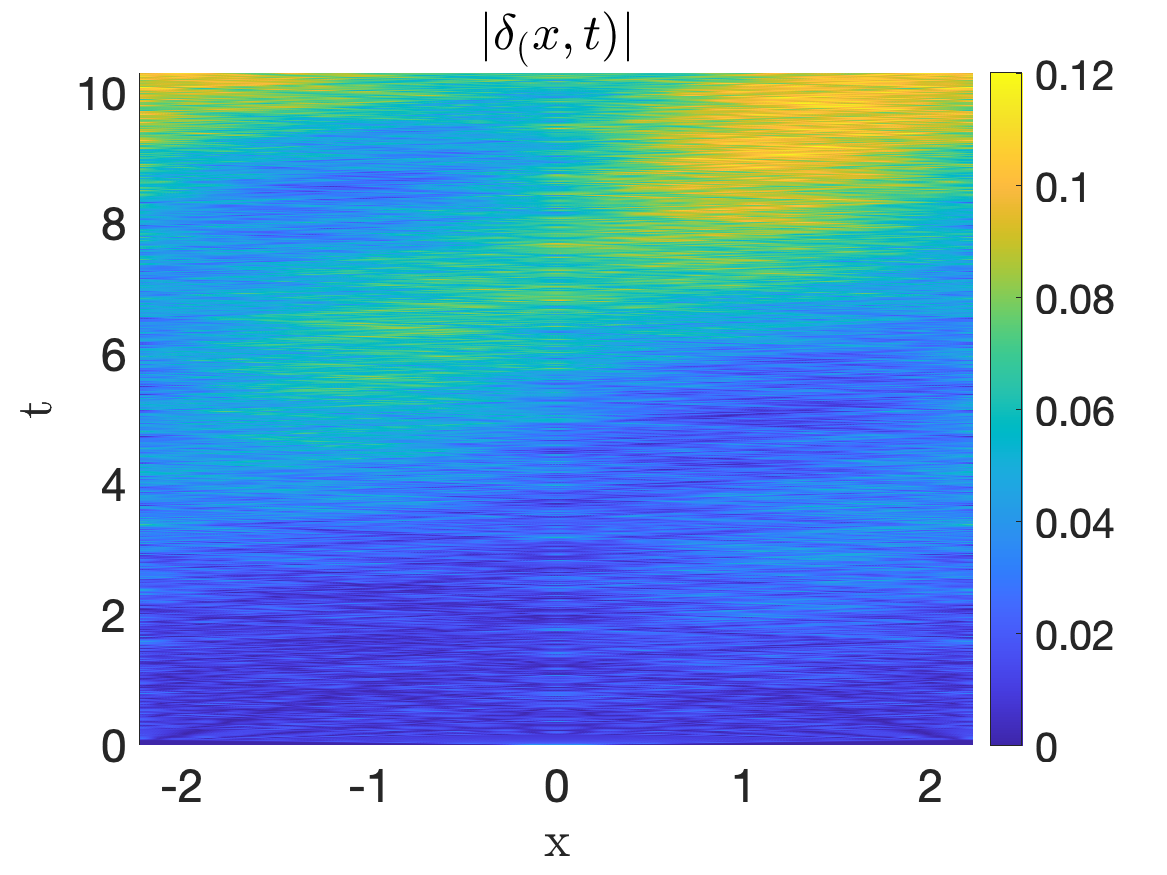} \,
\includegraphics[width=0.49\textwidth]{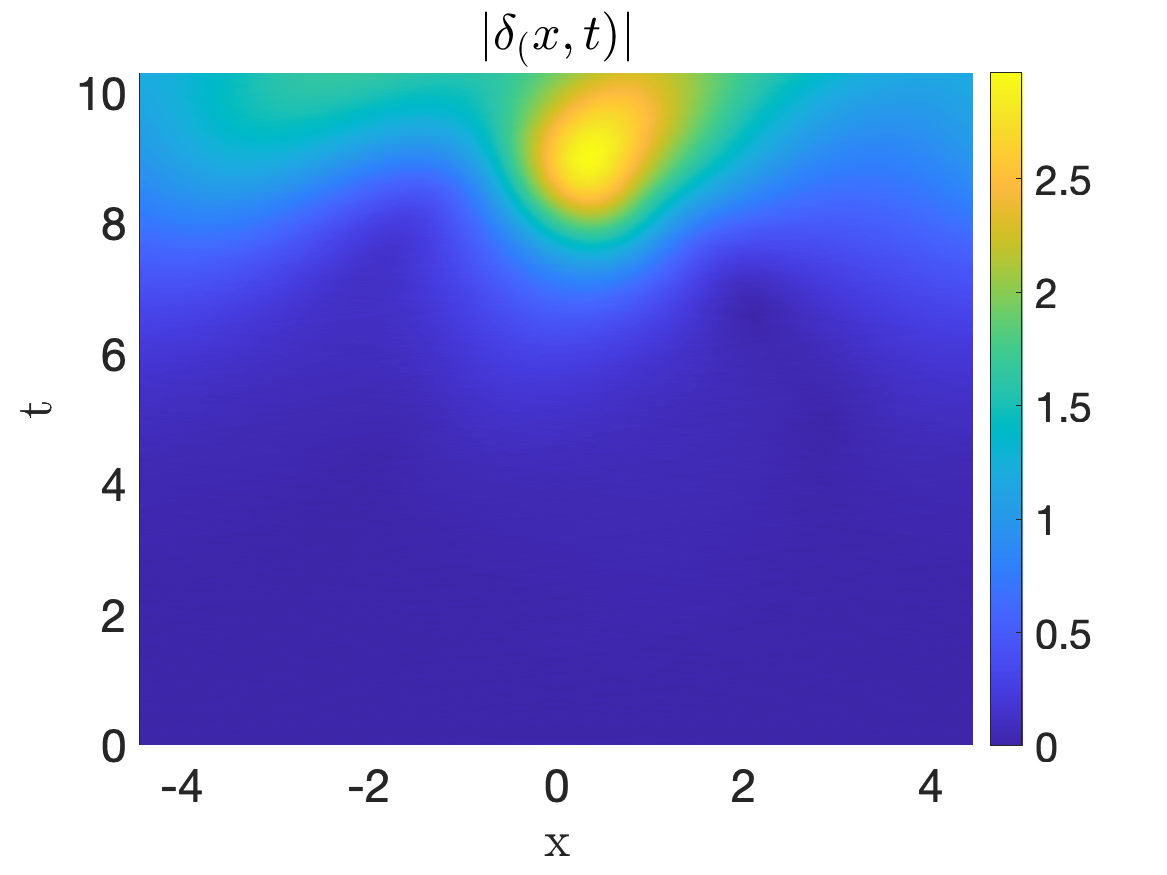} \\
\includegraphics[width=0.49\textwidth]{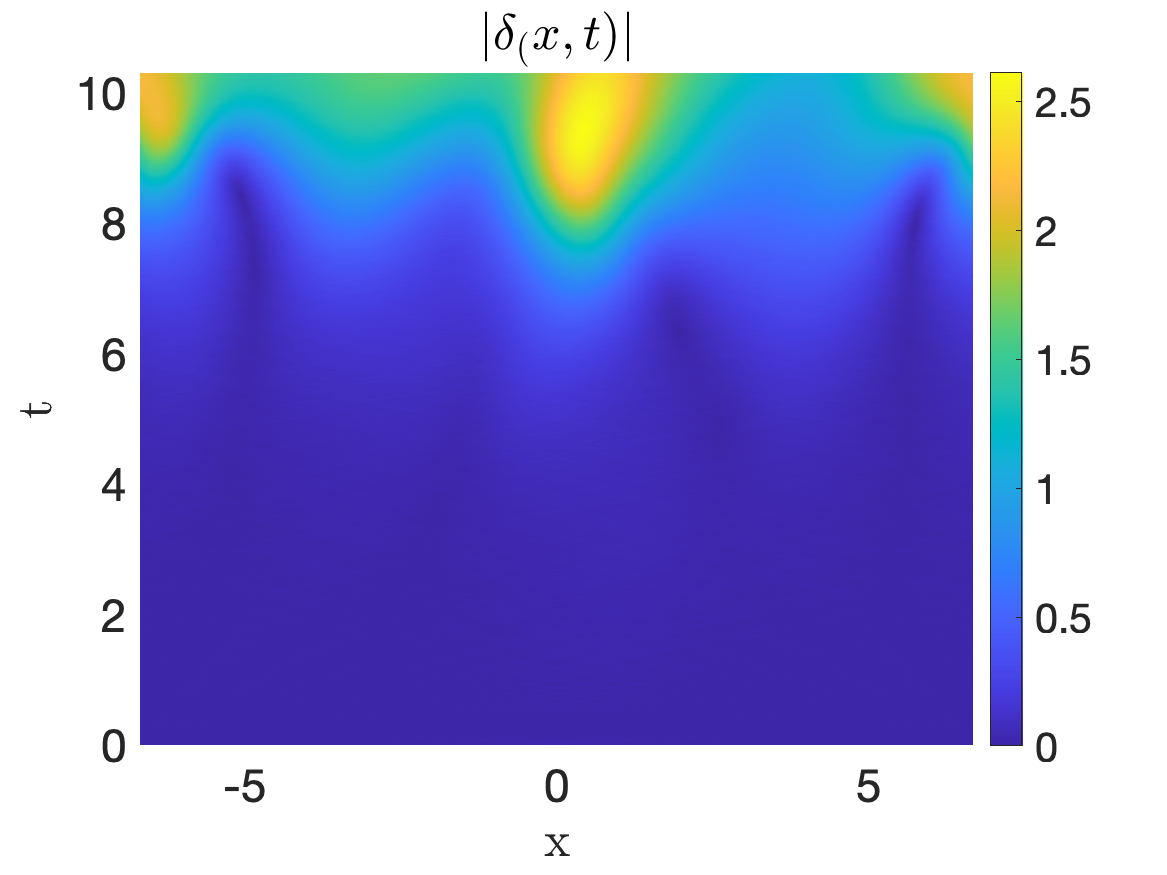} \,
\includegraphics[width=0.49\textwidth]{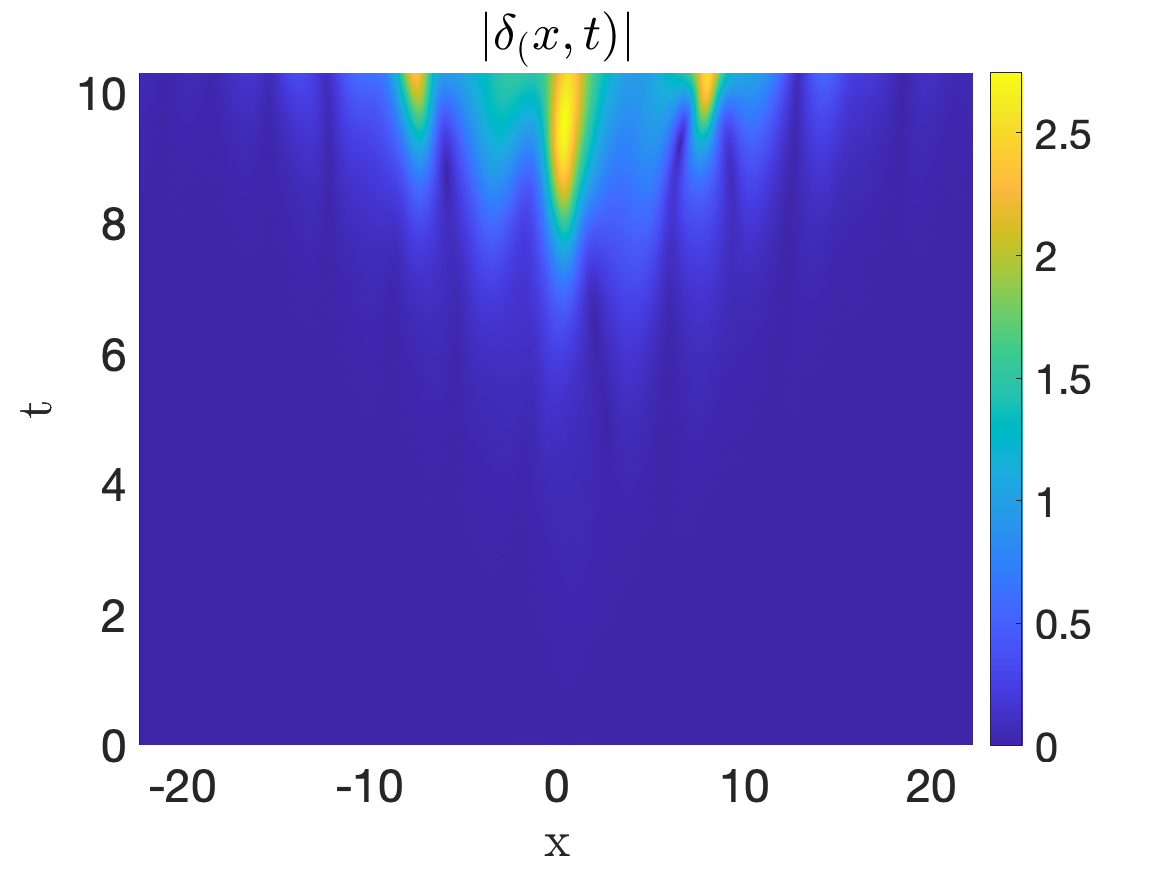}
\caption{Space-time plot of the modulus of the inhomogeneity, $|\delta(x,t)|,$ when propagated on intervals of different length $L$ according to equations \eqref{eq:87yu76yt6} -- \eqref{eq:87yu76yt62}. The reference length $L_0$ is reported in equation \eqref{eq:87yu76yt6}.  {\bf Top left:} $L=L_0.$ {\bf Top right:} $L=2L_0.$ {\bf Bottom left:} $L=3L_0.$ {\bf Bottom right:} $L=10L_0.$}
\label{fig:333}
\end{figure}

\subsection{Numerical investigation of the fully nonlinear MI} \label{sec:3}

In what follows we present numerical solutions of the problem
\begin{equation}\label{eq:nls4}
\begin{array}{c}
i \partial_t u + p \Delta u + q |u|^2u=0, \qquad u(x,0)=A(1+\delta_j(x))\\  u(-\frac{L}2,t)=u(\frac{L}2,t), \quad \partial_x u(-\frac{L}2,t)=\partial_x u(\frac{L}2,t),
\end{array}
\end{equation}
for the initial inhomogeneities 
\begin{equation}\label{eq:ICinh}
\begin{array}{c}
\delta_1(x) = 0.03 \cdot A \cdot \mathrm{sech}(15x) \cos(5x),\qquad   
\delta_2(x) = 0.03 \cdot A \cdot \mathrm{sech}(15x), \\
\delta_3(x) = 0.03 \cdot A \cdot e^{-3x^2},\quad 
\delta_4(x) = 0.03 \cdot A \cdot e^{-x^4}, \quad 
\delta_5(x) = 0.06 \cdot A \cdot x\cdot e^{-x^4}, 
\end{array}
\end{equation}
and $p=q=A=1.$  The bifurcation length in this case is $L_c\approx 4.45;$ the solutions are computed for $L=0.98 L_c,$ $L=1.3 L_c,$ $L=3L_c$ and $L=10 L_c.$  The computation time is $t\in [0,10].$ The numerical scheme used is a relaxation in time with second-order finite differences in space, and it satisfies mass and energy conservation on the discrete level \cite{Besse2004,Besse2021}. 
The inhomogeneity $\delta(x,t)$ for $t\in[0,10]$ is defined as 
\begin{equation}\label{eq:defddd}
\delta(x,t) := (u(x,t) - Ae^{i q A^2 t})\frac{1}Ae^{-i q A^2 t}
\end{equation}
consistently with equation \eqref{eq:4}. Numerical results can be found in Figure \ref{fig:b} 
as well as in Table \ref{tab:1}.

The initial inhomogeneities $\delta_1$ and $\delta_2$ are localized $\mathrm{sech}$ bumps, and are quite similar to each other. $\delta_3$ and $\delta_4$ are also bumps with different profiles, and $\delta_5$ is a localised wave. 
In Table \ref{tab:1} the maximum moduli of the inhomogeneities are recorded when propagated on intervals of different lengths. There is a clear confirmation of the abrupt bifurcation, namely an abrupt change in behaviour from $L=0.98L_c$ to $L=1.3L_c.$

On the larger intervals $L=3L_c,$ $L=10L_c,$ coherent structures emerge. These are a sign that the MI has been resolved in a manner qualitatively similar to what would happen on the real line (that is, until the structures reach close to the boundary of the computational domain). These structures all involve localised maxima supported on neighbourhoods of size roughly $L_c$ each. In the larger domain, the well-known space-time cone \cite{Biondini2016,Biondini2017,Biondini2018} is clearly visible. 

The presence of a space-time cone in the nonlinear phase of the MI, reported for a wide array of dispersive problems \cite{Biondini2016,Biondini2017,Biondini2018}, means that there is generically a finite speed of propagation of the boundary between the plane wave region (for $|x|\gg 1$) and the ``inhomogeneous'' region. Moreover, the ``inhomogeneous'' region seems to settle in a pattern determined by the equation and not by the specific initial condition (see also Figures 1 and A1 of \cite{Athanassoulis2023}). These patterns that appear within the cone seem to be related to so-called dispersive shock solutions \cite{Biondini2018a}.



Crucially, when the length of the domain becomes $L<L_c,$ we don't merely see a periodised or slightly off version of the cone, but something completely different:   the amplification is small or non-existent and there is no spatial pattern at all, see  left graph in Figure \ref{fig:b}. Indeed, for all initial data, when $L=0.98L_c,$ the inhomogeneity initially disperses rapidly and then slowly grows in homogeneous way over the whole computational domain. 


\subsection{An instance of gMI on an interval of length $L$}
\label{sec:nlgmit}

We set
\begin{equation}\label{eq:87yu76yt6}
	L_0=4.4518, \qquad  p=q=1
\end{equation}
and consider an initial condition 
\begin{equation}
	u_0(x) =  0.018  e^{-2\pi i \frac{x}{L_0}} +0.899  + 0.1252 e^{2\pi i \frac{x}{L_0}}.
\end{equation}
This  $u_0$ can be thought of as a perturbation of the plane wave solution \eqref{eqpwdef} at $t=0,$ and the lengthscale $L_0$ is numerically close to the $L_c$ of problem \eqref{eq:nls4} -- \eqref{eq:defddd}. The goal of this computation is to see whether a phenomenology analogous to what was seen in Section \ref{sec:3} arises beyond plane wave backgrounds. Thus the initial condition $u_0(x)$ plays the role of a realization of a very narrow spectrum.

In this setting the ``homogeneous background solution'' $u(x,t)$ is the solution of
\begin{equation}\label{eq:54u}
	\begin{array}{c}
i \partial_t u + p \Delta u + q |u|^2u=0, \qquad u(x,0)=u_0(x)\\  u(-\frac{L}2,t)=u(\frac{L}2,t), \quad \partial_x u(-\frac{L}2,t)=\partial_x u(\frac{L}2,t),
\end{array}
\end{equation}
and the inhomogeneously perturbed solution $V(x,t)$ is the solution of 
\begin{equation}\label{eq:54v}
	\begin{array}{c}
i \partial_t v + p \Delta v + q |v|^2v=0, \qquad v(x,0)=u_0(x) + 0.07 \mathrm{sech}(15x) \cos(5x)\\  v(-\frac{L}2,t)=v(\frac{L}2,t), \quad \partial_x v(-\frac{L}2,t)=\partial_x v(\frac{L}2,t).
\end{array}
\end{equation}
The inhomogeneity $\delta(x,t)$ is thus defined as
\begin{equation}\label{eq:87yu76yt62}
	\delta(x,t) = v(x,t) - u(x,t),
\end{equation}
and the main question is whether $\delta(x,t)$ grows significantly or not for different values of $L.$

Immediately when we consider background solutions other than plane waves, we observe that the length $L$ and the initial conditions $u_0,$ $\delta$ have to be compatible: the homogeneous $u_0$ has to satisfy periodic BCs on $L,$ hence in this case we will consider only integer multiples of $L_0.$ On the other hand, the localized $\delta(x,0)$ has to be numerically zero at the ends of the domain, up to satisfactory accuracy. By taking $L=NL_0,$ for $N\in\{1,2,3,10\}$ both of these requirements are satisfied. The final time is taken to be $T=10.3.$

The numerical scheme used is the same as in Section \ref{sec:3} (a relaxation in time with second-order finite differences in space, and it satisfies mass and energy conservation on the discrete level \cite{Besse2004,Besse2021}). The time step used is $dt=4 \cdot 10^{-3}$ and the space mesh size is $dx=4\cdot 10^{-3}.$ The modulus $|\delta(x,t)|$ can be found in Figure \ref{fig:333}, and it is clear that for the short enough computational domain $L=L_0$  the growth of the inhomogeneity is  completely suppressed (in particular the inhomogeneity in the $L^\infty$ norm is not even doubled). On the other hand, the space-time cone, inside which the inhomogeneity grows by more than an order of magnitude is very clearly visible when $L=10L_0.$

\begin{figure}
	\includegraphics[width=0.49\textwidth]{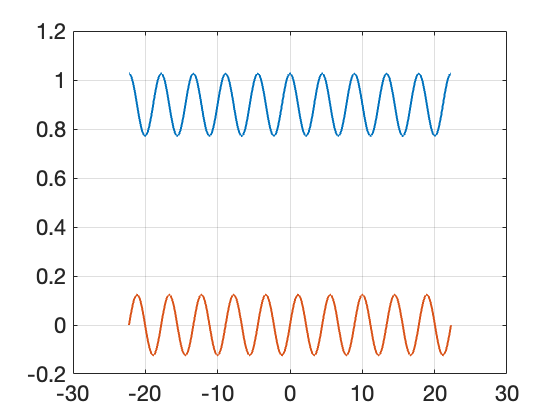}\,
		\includegraphics[width=0.49\textwidth]{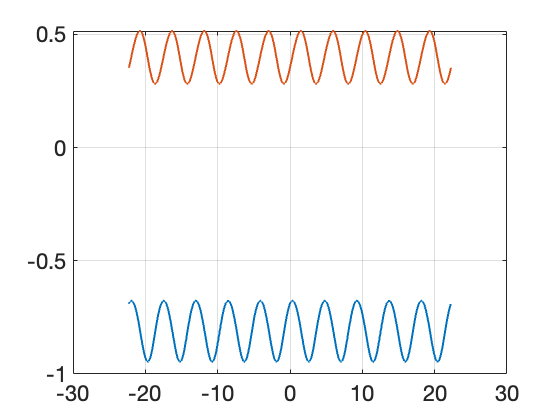}
		\caption{Plots of the real part (blue) and imaginary part (red) of the background solution $u,$ as in equation \eqref{eq:54u}, on the interval with length $L=10L_0.$ {\bf Left:} for $t=0.$ {\bf Right:} for $t=T.$}
		\label{fig:44}
\end{figure}

\begin{figure}
	\includegraphics[width=0.49\textwidth]{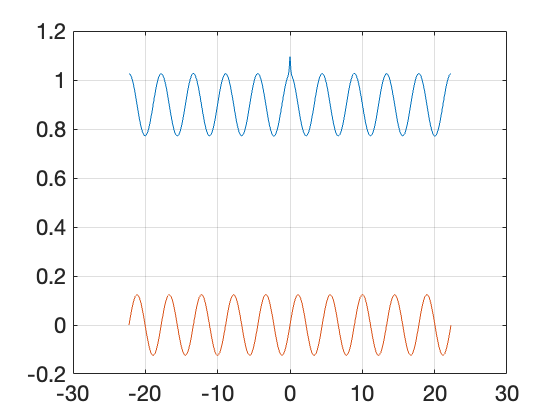}\,
		\includegraphics[width=0.49\textwidth]{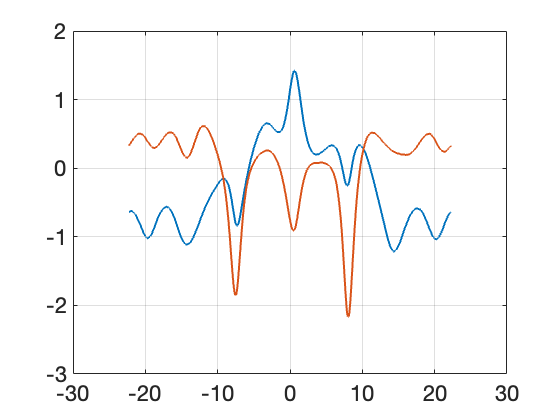}
		\caption{Plots of the real part (blue) and imaginary part (red) of the perturbed solution $v,$ as in equation \eqref{eq:54v}, on the interval with length $L=10L_0.$ {\bf Left:} for $t=0.$ {\bf Right:} for $t=T.$}
		\label{fig:55}
\end{figure}

 To provide some more context, the background solution and the inhomogeneity are plotted in the intitial and final time in Figures \ref{fig:44} and \ref{fig:55}. Observe how the background solution qualitatively behaves similarly to the plane wave solution.

\section{Conclusions} \label{sec:conc}

In this paper we presented for the first time a construction of the random phase approximation in the NLS equation that comes with a convergence result: using the scaling of equation \eqref{eq:scaling}, for large enough $L,$ the instability condition of the truncated problem \eqref{eq:nls2a} converges to that of the original problem \eqref{eq:nls1}. 

In the proof it is clear that equation \eqref{eq:scaling} is not the only way to prepare initial data that would have this property.
There are other ways to prepare initial data in the literature, e.g. with randomized amplitudes as well as phases of the modes \cite{Athanassoulis2021}, and it is possible that these could also lead to a consistent stability condition. On the other hand, recent  results in different regimes  \cite{el2023,Buckmaster2021,ZHani23}  highlight that different ways to prepare the  initial conditions may lead to genuinely different outcomes. 

The importance of using large enough computational domains is well understood for the classical MI. This was further elaborated in Section \ref{sec:2}, where a lengthscale of $L_c=O(\lambda_0^2)$ was found to be necessary for the MI to be present even for plane waves. If $\lambda_0$ is in the hundreds of meters, as in ocean waves, this could mean $L_c$ hundreds of wavelengths long.

The role of $L$ in the gMI is less well understood. In works where the onset of gMI was actively investigated, it has been reported that long computational domains are needed; indeed in \cite{Ribal2013a},  it was reported that a computational domain of  $15L_c$ was needed for a JONSWAP sea state (where $L_c$ is the corresponding critical length for the plane wave with the same wavelength and wave height as the peak values for the sea state). Around $130$ wavelengths were found to be required in \cite{Athanassoulis2023}, while in \cite{Onorato2003} the computational domain is about $100$ wavelengths.
On the other hand, in many papers the question of the length of computational domain seems to not be addressed in any detail. While more work is needed to conclusively resolve this question for realistic ocean wave spectra, we have demonstrated here (in Section \ref{sec:3} and Figure \ref{fig:333}) that the capacity of short computational domains to suppress the growth of inhomogeneities extends beyond plane wave backgrounds.


This work was carried out assuming NLS dynamics. This is of course an approximation for ocean waves, 
but it should be noted that recent rigorous works have shown that   the MI for the fully nonlinear water waves problem is accurately captured by the nonlinear Schr\"odinger equation \cite{berti22b,strauss23}. Furthermore, the analysis carried out here for the truncated NLS problem can be carried out for the truncated Zakharov problem in the spirit of \cite{Andrade2020b}, thus producing a stability relation for the periodized CSY equation (which is 2-dimensional and broadband, since it assumes  Zakharov dyamics).
This would also allow the investigation of crossing seas scenarios, either in fully 2D through the aforementioned CSY equation, or using an Alber system for the interaction of two unidirectional wavetrains at an angle \cite{Gramstad2018a,Athanassoulis2021}. Crossing seas are very important with regard to real-world marine safety, and much more can be done to investigate them \cite{Davison2022,Onorato2010,Steer2019}. More broadly, incorporating more physics into the model, such as wind input \cite{Athanassoulis2017,Toffoli2023} and vorticity \cite{Curtis2020} seems to be a very promising direction with regard to real-world ocean waves.

\medskip
\noindent {\bf Acknowledgment.}  AA would like to thank Prof. G. A. Athanassoulis, Prof. T. Sapsis, Dr. O. Gramstad and Dr. T. Tang for helpful discussions.

\smallskip

\bibliographystyle{siam}
\bibliography{MyCollection.bib}

\appendix

\end{document}